\documentclass[11pt]{amsart}

\usepackage{amsmath,amssymb,geometry}
\usepackage{verbatim}
\usepackage{pxfonts}
\usepackage{hyperref}
\hypersetup{colorlinks,citecolor=blue,plainpages=false,hypertexnames=false}
\newtheorem{theorem}{Theorem}[section]
\newtheorem{proposition}[theorem]{Proposition}
\newtheorem{lemma}[theorem]{Lemma}

\theoremstyle{definition}
\newtheorem{definition}[theorem]{Definition}

\theoremstyle{remark} 
\newtheorem{remark}[theorem]{Remark}

\numberwithin{equation}{section}

\begin{document}

\title[On an asymptotic characterisation of Griffiths semipositivity]{On an asymptotic characterisation of Griffiths semipositivity}

\author{Apoorva Khare}
\address{Department of Mathematics, Indian Institute of Science; and
Analysis \& Probability Research Group; Bangalore -- 560012, India}
\email{khare@iisc.ac.in}

\author{Vamsi Pritham Pingali}
\address{Department of Mathematics, Indian Institute of Science; and
Analysis \& Probability Research Group; Bangalore -- 560012, India}
\email{vamsipingali@iisc.ac.in}

\begin{abstract} 
We prove that certain possibly non-smooth Hermitian metrics are
Griffiths-semipositively curved if and only if they satisfy an asymptotic
extension property. This result answers a question of
Deng--Ning--Wang--Zhou in the affirmative.
\end{abstract}
\maketitle
\section{Introduction}\label{Introsec}

The positivity of the curvature $\Theta$ of the Chern connection of a
Hermitian holomorphic vector bundle $(E,h)$ over a complex manifold $M$
plays an important role in algebraic geometry through extension problems.
For a line bundle, there is only one notion of positivity, namely, the
curvature being a K\"ahler form. For a vector bundle, there are several
competing inequivalent notions, of which the most natural are
\emph{Griffiths positivity} ($\langle v, \sqrt{-1}\Theta v\rangle$ is a
K\"ahler form) and \emph{Nakano positivity} (the bilinear form defined by
$\sqrt{-1}\Theta$ on $T^{1,0}M \otimes E$ is positive-definite). A famous
conjecture of Griffiths \cite{Griffiths} asks whether \emph{ample vector
bundles} ($\mathcal{O}_E(1)$ over $\mathbb{P}(E)$ admits a
positively-curved metric) admit Griffiths-positively curved metrics. The
conjecture is still open. However, a considerable amount of work has been
done to provide evidence in its favour \cite{Bo, Camp, dem2, demskoda,
Diverio, Guler, MT, Naumann, Pinchern, Um}.

In the K\"ahler case, Demailly--Paun \cite{dp} proved a
Nakai--Moizeshon--type criterion to characterise K\"ahler classes.
Despite the assumptions and conclusions involving smooth objects, their
proof used singular objects like positive currents and singular K\"ahler
potentials. A similar phenomenon might play a role in the study of the
Griffiths conjecture and hence it is fruitful to study singular Hermitian
metrics on vector bundles. This topic has also been well-studied
\cite{decat, Rauf, BernPaun,Bernd, Hosono,Pauntakayama,Ra2,RaufiChern}
and seems to hold some surprises. For instance, even if a bundle is
Griffiths-positively curved (in a certain sense), the curvature may not
exist as a current (Theorem $1.3$ in \cite{Rauf}).

In the quest for alternate characterisations of these notions of
positivity, and defining similar notions for singular metrics, the
following definition \cite{Deng3} involving the asymptotics  of
$L^2$-extension constants has proved to be a useful measure of positivity
\cite{Deng1, Deng2, Deng3, Hosin, Ini}.

\begin{definition}
Let $E$ be a holomorphic vector bundle over an $n$-dimensional complex
manifold $X$. A singular Hermitian metric $h$ is said to satisfy the
multiple coarse $L^2$-extension property if the following hold.
\begin{enumerate}
\item For every open subset $D\subset X$ and every holomorphic section
$s: D\rightarrow E^*\vert_D$ that is not identically zero, the function
$\ln \Vert s \Vert_{h^*}^2$ is upper-semicontinuous.

\item Consider any cover of $M$ by relatively compact Stein trivialising
coordinate neighbourhoods of the form $(\Omega'' \subset M, z, \{e_i\})$
and a subcover of Stein open subsets $\Omega' \Subset \Omega \Subset
\Omega''$. Then, for every integer $m\geq 1$ there exists a constant
$C_m$  satisfying the following conditions:
\begin{enumerate}
\item Subexponential growth: $\displaystyle \lim_{m \to \infty} \frac{\ln
C_m}{m} = 0$. 
\item Controlled extension: If $p\in \Omega'$ and $a\in E_p$ with $\Vert
a \Vert_h <\infty$, for every integer $m\geq 1$, there exists a
holomorphic extension $f_m : \Omega \rightarrow E^{\otimes m}$ of
$a^{\otimes m}$ (that is, $f_m (p)=a^{\otimes m}$) satisfying 
\begin{equation}\label{eq:l2extensionproperty}
 \int_{\Omega'} \Vert f_m \Vert_{h^{\otimes m}}^2
 \frac{(\sqrt{-1}\partial \bar{\partial} \vert z \vert^2)^n}{n!} \leq C_m
 \Vert a \Vert_{h(p)}^{2m}.
\end{equation}
\end{enumerate}
\end{enumerate}
\label{def:main}
\end{definition}

We note that in the definition above, $a^{\otimes m} \in S^m E_p$. It
turns out (Lemma \ref{lem:fixedpoint}) that the metric on $S^m E$ induced
from $E$ is actually the same as the metric induced by an orthogonal
projection from $E^{\otimes m}$. This observation motivates the following
generalisation of Definition \ref{def:main} : A singular Hermitian metric
$h$ is said to satisfy the \emph{general} multiple coarse $L^2$-extension
property if it satisfies Definition \ref{def:main} with $f_m(p)$ being
\emph{any} given element of $S^m E_p$, i.e., if $b_m \in S^m E_p$ then
there exists a holomorphic section $f_m :\Omega\rightarrow S^m E$ such
that $f_m(p)=b_m$ and 
\begin{equation}\label{eq:l2extensionpropertygeneral}
 \int_{\Omega'} \Vert f_m \Vert_{h^{\otimes m}}^2
 \frac{(\sqrt{-1}\partial \bar{\partial} \vert z \vert^2)^n}{n!} \leq C_m
 \Vert b_m \Vert_{S^m h(p)}^2 \ ,
\end{equation}
where $C_m$ satisfies subexponential growth.

For the remainder of this paper, unless specified otherwise, an integral
over a coordinate chart is understood to be an integral with the
Euclidean volume form, similar to \eqref{eq:l2extensionproperty}.

\begin{remark}
The definition given in \cite{Deng3} differs slightly from Definition
\ref{def:main} in two (minor) aspects. Firstly, the definition in
\cite{Deng3} is for general Finsler metrics. Secondly, the controlled
extension property in our definition requires control over the extension
(to a set $\Omega$) on a smaller set $\Omega'$.
\label{rem:afterdef}
\end{remark}

We recall that a singular Hermitian metric $h$ is said to be
Griffiths-semipositively curved if whenever $u$ is a local holomorphic
section of $E^*$, the function $\vert u \vert^2_{h^*}$ is a
plurisubharmonic (psh) function \cite{Rauf}. (It turns out that this
definition is equivalent to $\ln \vert u \vert^2_{h^*}$ being psh.) In
\cite{Deng1} (Theorem 6.4), it was proved that multiple coarse
$L^2$-extension implies Griffiths semipositivity. Since the proof is
local, it is easily seen to apply even to our definition. A question was
raised as to whether it completely characterises Griffiths
semipositivity. We answer that question in the affirmative. Slightly more
strongly:

\begin{theorem}
Let $h$ be a singular Griffiths semipositively curved Hermitian metric on
a holomorphic vector bundle $E$ over a complex Hermitian manifold
$(X,\omega)$. Let $h_0$ be a fixed smooth background metric on $E$. If
$\ln\det(hh_0^{-1})$ is bounded on compact sets, then $(E,h)$ satisfies
the general multiple coarse $L^2$-extension property.
Moreover, one can choose a uniform extension constant $C = C_m$ that is
independent of $m$.
\label{thm:main}
\end{theorem}

\begin{remark}
If $h$ is continuous, it trivially meets the requirements of Theorem
\ref{thm:main}. If instead, $0\leq \sqrt{-1}\bar{\partial}\partial \ln
(\det(h)) \leq C\omega$, then $\ln(\det(h))$ and $-\ln(\det(h))$ are
quasi-psh and hence satisfy the hypotheses of Theorem \ref{thm:main}. 
\end{remark}

\begin{remark}
The main point of Theorem \ref{thm:main} is that, while typically Nakano
semipositivity produces extension theorems, this theorem merely needs
Griffiths semipositivity.
\end{remark}

It is interesting to know whether our result can be improved to general
singular Hermitian metrics. However, seeing that $\det(h)$ seems to play
an important role, we are pessimistic about such a result. There are
other positivity notions for vector bundles like MA-positivity, for
instance \cite{Pin}. It might be fruitful to explore a similar
extension/estimate-type characterisation for such notions as well.

\subsection*{Acknowledgements}
This work is partially supported by
grant F.510/25/CAS-II/2018(SAP-I) from UGC (Govt. of India),
Ramanujan Fellowship grant SB/S2/RJN-121/2017,
MATRICS grants MTR/2017/000295 and MTR/2020/000100, 
and SwarnaJayanti Fellowship grants SB/SJF/2019-20/14 and
DST/SJF/MS/2019/3 from SERB and DST (Govt. of India). We are grateful to
the anonymous referee for useful suggestions to improve the paper.

\section{Proof}\label{Proof}

We first prove Theorem \ref{thm:main} for smooth Hermitian metrics.

\begin{proposition}\label{prop:2implies1insmooth}
Let $h_0$ be a smooth background metric. If $h$ is a smooth Hermitian
metric, then Griffiths semipositivity implies the general multiple coarse
extension property.  Moreover, if $h_{\nu}$ is a family of such smooth
Hermitian metrics and there exists a constant $K$ such that
$\frac{1}{K}\det(h_0) \leq \det(h_{\nu})\leq K\det(h_0)$, then the
extension constants $C_m$ depend only on $K, h_0,$ and the chosen
coordinate trivalising charts.
\end{proposition}

\begin{proof}
Consider a point $p\in \Omega' \Subset \Omega \Subset \Omega'' $ and $b_m
\in S^m E_p$. We want to extend $b_m$ to a holomorphic section in
$\Omega$ with $L^2$-estimates in $\Omega'$. Let $\epsilon_m>0$ be a
sequence of real numbers. The new metric
$\tilde{h}_{\epsilon_m}:=he^{-\epsilon_m \vert z \vert^2}$ is strictly
Griffiths-positive on $\Omega$.

We recall that the symmetric group $S_m$ acts on $E^{\otimes m}$ which
decomposes into irreducible representations. The metric
$\tilde{h}_{\epsilon_m}$ induces an $S_m$-invariant metric on $E^{\otimes
m}$. The fixed point set is $S^m E$ and hence we decompose $E^{\otimes m}
= S^m E \oplus V$. Note that $V$ is stable under the action of $S_m$.

\begin{lemma}
The subbundle $S^m E$ is orthogonal in the induced metric to $V$.
\label{lem:fixedpoint}
\end{lemma}

\begin{proof}
Indeed, suppose $x\in S^m E_q$, $y \in V_q$ (for some $q$). Then
$m!\langle x,y \rangle =  \sum_{g\in S_m} \langle g \cdot x ,y \rangle =
\sum_{g \in S_m} \langle x, g^{-1} \cdot y \rangle = \langle x, y_0
\rangle$ where $y_0 = \sum_{g \in S_m} g^{-1} \cdot y$ is in $V_q$ as
well as in the fixed-point set $S^m E_q$. Hence, $y_0 =0$ and so is
$\langle x,y\rangle$.
\end{proof}   

Endow $\Omega$ with the Euclidean metric. Since $E$ is trivial over
$\Omega$, we pretend that $\det(E)$ is a trivial bundle. Let $r$ be the
rank of $E$. At this point, suppose $h_{\nu}$ is a family of smooth
Hermitian metrics as in the statement of the proposition, and let
$\epsilon_m=\frac{1}{m}$. We drop the subscript $\nu$ for the remainder
of the proof.

A result of Demailly--Skoda \cite{demskoda} shows that if $(E,h)$ is
Griffiths-positively curved, then $E \otimes \det(E)$ with the induced
metric is Nakano-positively curved. This result was generalised in
Theorem $7.2$ of \cite{Liu} which states that the induced metric on $S^m
E\otimes \det(E)$ is Nakano-positively curved for all $m\geq 1$.

By the Ohsawa--Takegoshi theorem for vector bundles, there exists a
universal constant $C$ (whose optimal value can be computed \cite{Guan})
and an extension $f_{m}$ of $b_m$ such that
\begin{equation}\label{eq:OTkey}
\int_{\Omega} \Vert f_{m} \Vert_{S^m
\tilde{h}_{\epsilon_m}}^2 \det(\tilde{h}_{\epsilon_m}) \leq C \Vert
b_m\Vert_{S^m \tilde{h}_{\epsilon_m}(p)}^2
\det(\tilde{h}_{\epsilon_m}).
\end{equation}

By Lemma \ref{lem:fixedpoint}, $\Vert b \Vert_{S^m h} = \Vert b
\Vert_{h^{\otimes m}}$ if $b\in S^m E$. Note that 
\[
\Vert f_m \Vert^2_{S^m \tilde{h}_{\epsilon_m}}
\det(\tilde{h}_{\epsilon_m}) = \Vert f_m \Vert_{S^m h}^2 e^{-m\epsilon_m
\vert z \vert^2} \det(h) e^{-r\epsilon_m \vert z \vert^2}
= \Vert f_m \Vert_{S^m h}^2 e^{-\left(1+\frac{r}{m}\right) \vert z
\vert^2} \det(h).
\]
Rewriting \eqref{eq:OTkey}, 
\begin{align*}
\frac{\inf_{\overline{\Omega}} \det(h_0) e^{-(r+1) \vert z \vert^2}}{K}
\int_{\Omega}\Vert f_{m} \Vert_{S^m h}^2 \leq &\ \int_{\Omega} \Vert f_m
\Vert_{S^m h}^2 e^{-\left(1+\frac{r}{m}\right) \vert z \vert^2} \det(h)\\
\leq &\ C \Vert b_m \Vert_{S^m h(p)}^2 e^{-\left(1+\frac{r}{m}\right)
\vert z(p) \vert^2} K \sup_{\overline{\Omega}} \det(h_0).
\end{align*}
Hence, we are done.
\end{proof}

Now we prove Proposition \ref{prop:2implies1insmooth} in greater
generality. Let $p \in \Omega'$ and $b_m\in S^m E_p$. Proposition 6.2 of
\cite{Rauf} implies that the duals of the convolutions of the dual metric
with mollifiers $\rho_{\nu}$, i.e., $h_{\nu}=((h^*)*\rho_{\nu})^*$ of $h$
(where $0< \nu\leq 1$), increase to $h$ pointwise as $\nu\rightarrow 0^+$
and are Griffiths-semipositively curved. Choose $\nu \leq \nu_0$ to be
small enough that the convolutions are defined on $\Omega_{2 \delta}$,
where $\Omega_{2 \delta} \Subset \Omega''$ is a $2 \delta$-neighbourhood
of $\Omega$ for some fixed small $\delta > 0$. By Theorem 7.2 of
\cite{Liu}, $S^m E\otimes \det(E)$ equipped with the metric induced from
$h_{\nu}$ is Nakano-semipositively curved.

Since $0<\nu\leq \nu_0$, by the monotonicity of $h_{\nu}$, we see that
$\frac{1}{L}\ln(\det(h_0)) \leq \ln (\det (h_{\nu_0})) \leq
\ln(\det(h_{\nu}))\leq \ln (\det(h)) \leq L \ln(\det(h_0))$ for some
constant $L$ independent of $\nu$. Hence we may use Proposition
\ref{prop:2implies1insmooth} to conclude that there exist extensions
$f_{m,\nu}$ of $b_m$ on $\Omega_{2 \delta}$ such that
\[
\int_{\Omega_{\delta}} \Vert f_{m,\nu}
\Vert_{h_{\nu_0}^{\otimes m}}^2 \leq \int_{\Omega_{\delta}} \Vert f_{m,\nu}
\Vert_{h_{\nu}^{\otimes m}}^2   \leq  C \Vert b_m
\Vert_{h_{\nu}^{\otimes m}(p)}^2 \leq C \Vert b_m
\Vert_{h^{\otimes m}(p)}^2
\]
where $C$ is a constant independent of $m$ and $\nu$. Henceforth, all
such constants will be denoted by $C$. Therefore, $f_{m,\nu}$ is
uniformly bounded (independent of $\nu$) in $L^2(\Omega_{\delta})$. The
sub-mean value property shows that it is pointwise bounded in
$\Omega_{\delta/2}$. Cauchy's estimates show that $\Vert f_{m,\nu}
\Vert_{C^3(\Omega_{\delta/3})}$ is bounded above by some $K_m$ (uniformly
in $\nu$). Thus, by the Arzela--Ascoli theorem, there is a sequence
$\nu_i \rightarrow 0$ such that $f_{m,\nu_i} \rightarrow v_m$ in
$C^2(\Omega_{\delta/4})$. The limit $v_m$ is a holomorphic extension of
$b_m$ over $\Omega$. Fixing $\nu_0$, by uniform convergence of
$f_{m,\nu_i}$, we see that
\begin{gather}
\int_{\Omega'} \Vert v_m \Vert_{h_{\nu_0}^{\otimes m}}^2 \leq C
\Vert b_m \Vert_{h^{\otimes m}(p)}^2. \nonumber 
\end{gather} 
By the monotone convergence theorem,
\begin{gather}
\int_{\Omega'} \Vert v_{m} \Vert_{h^{\otimes m}}^2 \leq C \Vert
b_m \Vert_{h^{\otimes m}(p)}^2. \nonumber
\end{gather}
Thus, the general multiple coarse $L^2$-extension property is met and the proof
of Theorem \ref{thm:main} is complete. \qed

\end{document}